\documentclass[11pt]{article}
\usepackage{amsthm, amsmath, amssymb, amsfonts, url, booktabs, tikz, setspace, fancyhdr, bm}
\usepackage{geometry}
\usepackage{hyperref, enumerate}
\usepackage[shortlabels]{enumitem}
\usepackage[babel]{microtype}
\usepackage[english]{babel}
\usepackage[capitalise]{cleveref}
\usepackage{comment}
\usepackage{bbm}
\usepackage{csquotes}
\usepackage{mathabx}
\usepackage{tikz}
\usepackage{graphicx}
\usepackage{float}
\usepackage{xcolor}
\usetikzlibrary{positioning, arrows.meta, shapes.geometric}

\counterwithin{figure}{section}

\newtheorem{theorem}{Theorem}[section]

\newtheorem{lemma}[theorem]{Lemma}

\newtheorem{claim}[theorem]{Claim}

\usetikzlibrary{decorations.pathmorphing}

\theoremstyle{definition}

\newtheorem*{defn-non}{Definition}

\newlist{Case}{enumerate}{2}
\setlist[Case, 1]{%
    label           =   {\bfseries Case \arabic*.},
    labelindent=1em ,labelwidth=1.3cm, labelsep*=1em, leftmargin =!
}
\setlist[Case, 2]{%
    label           =   {\bfseries Subcase \arabic{Casei}.\arabic*.},
    labelindent=-1em ,labelwidth=1.3cm, labelsep*=1em, leftmargin =!
}

\newenvironment{poc}{\begin{proof}[Proof of claim]}{\end{proof}}

\usepackage{todonotes}

\title{Largest dyadic dual VC-dimension of non-piercing families} 
\author{
Xinqi Huang\thanks{School of Mathematical Sciences, University of Science and Technology of China, Hefei, China and Extremal Combinatorics and Probability Group (ECOPRO), Institute for Basic Science (IBS), Daejeon, South Korea.
Email:huangxq@mail.ustc.edu.cn.
}
\and
Yuzhen Qi\thanks{School of Mathematics, Shangdong University, Jinan, China, and Extremal Combinatorics and Probability Group (ECOPRO), Institute for Basic Science (IBS), Daejeon, South Korea. Email: yzq\_sdu\_edu@163.com. }
\and
Mingyuan Rong\thanks{School of Mathematical Sciences, University of Science and Technology of China, Hefei,
China.
Email: rong\_ming\_yuan@mail.ustc.edu.cn.}
\and
Zixiang Xu\thanks{Extremal Combinatorics and Probability Group (ECOPRO), Institute for Basic Science (IBS), Daejeon, South Korea. Email: zixiangxu@ibs.re.kr.}
}
\begin{document}
\date{}
\maketitle
\begin{abstract}
The dyadic dual VC-dimension of a set system \( \mathcal{F} \) is the largest integer \( \ell \) such that there exist \( \ell \) sets \( F_1, F_{2}, \dots, F_\ell \in \mathcal{F} \), where every pair \( \{i, j\} \in \binom{[\ell]}{2} \) is witnessed by an element \( a_{i,j} \in F_i \cap F_j \) that does not belong to any other set \( F_k \) with \( k \in [\ell] \setminus \{i, j\} \). In this paper, we determine the largest dyadic dual VC-dimension of a non-piercing family is exactly $4$, providing a rare example where the maximum of this parameter can be determined for a natural family arising from geometry. As an application, we give a short and direct proof that the transversal number \( \tau(\mathcal{F}) \) of any non-piercing family is at most \(C\nu(\mathcal{F})^9 \), where \( \nu(\mathcal{F}) \) is the matching number and $C$ is a constant. This improves a recent result of P\'{a}lv\"{o}lgyi and Z\'{o}lomy.
\end{abstract}

\section{Introduction}
Given a ground set \( V \) and a set system \( \mathcal{F} \subseteq 2^V \), the \emph{transversal number} \( \tau(\mathcal{F}) \) is the smallest integer \( t \) such that there exists a subset \( T \subseteq V \) of size \( t \) intersecting every set in \( \mathcal{F} \). The \emph{matching number} \( \nu(\mathcal{F}) \) is the maximum number of pairwise disjoint sets in \( \mathcal{F} \). It is clear that \( \nu(\mathcal{F}) \le \tau(\mathcal{F}) \), but in general, \( \tau(\mathcal{F}) \) cannot be bounded from above by any function of \( \nu(\mathcal{F}) \). Thus, it is of particular interest to identify specific set systems for which such an upper bound does exist.

A \emph{region} is a connected, compact subset of the plane whose boundary consists of finitely many disjoint Jordan curves. A family of regions is said to be in \emph{general position} if the boundaries of any two regions intersect in only finitely many points. A family \( \mathcal{F} \) of regions is called \emph{non-piercing} if, for any two regions \( F, G \in \mathcal{F} \), the difference \( F \setminus G \) is connected. Intuitively, this means that removing one region from another does not disconnect it. In a recent work, P\'{a}lv\"{o}lgyi and Z\'{o}lomy~\cite{2025Nonpierc} showed that for any family \( \mathcal{F} \) of non-piercing regions, there exists a function \( f \) such that \( \tau(\mathcal{F}) \le f(\nu(\mathcal{F})) \). Their proof cleverly combines a variant of the \((p, q)\)-theorem due to Matou\v{s}ek~\cite{2004DCGMatousek}, VC-dimension arguments, and the Clarkson--Shor probabilistic method~\cite{1989DCG}.

In contrast, we present a more direct approach based on a classical result of Ding, Seymour, and Winkler~\cite{1994DingSeymour}. This method yields a polynomial bound on the transversal number without invoking the full machinery of VC-dimension theory. More precisely, we obtain the following result.

\begin{theorem}\label{thm:main}
If \( \mathcal{F} \) is a non-piercing family of regions, then \( \tau(\mathcal{F}) \le C\nu(\mathcal{F})^{9} \) for some constant $C>0$.
\end{theorem}
As a direct consequence, the disjointness graph of any non-piercing family is polynomially \( \chi \)-bounded, and in particular, satisfies the Erd\H{o}s--Hajnal property. P\'{a}lv\"{o}lgyi and Z\'{o}lomy~\cite{2025Nonpierc} also applied this result to upper bound the maximum length of a facial cycle in a planar graph of girth at least \( \ell \) that is inclusion-maximal. For further details on this topic, we refer the readers to~\cite{2025Axen,2025Nonpierc}.

We now introduce the central parameter of this paper. For a set system \( \mathcal{F} \), the \emph{dyadic dual VC-dimension} \( \lambda(\mathcal{F}) \) is the largest integer \( \ell \) such that there exist sets \( F_{1}, F_{2}, \ldots, F_{\ell} \in \mathcal{F} \) satisfying the following property: for every pair \( \{i, j\} \in \binom{[\ell]}{2} \), there exists an element \( a_{i,j} \in F_i \cap F_j \) such that \( a_{i,j} \notin F_k \) for any \( k \in [\ell] \setminus \{i, j\} \). In 1994, Ding, Seymour, and Winkler~\cite{1994DingSeymour} established a connection between the transversal number and matching number via the dyadic dual VC-dimension.
\begin{lemma}[\cite{1994DingSeymour}]\label{1994}
For any set system \( \mathcal{F} \), we have
\[
\tau(\mathcal{F}) \le 11\lambda(\mathcal{F})^2 \cdot (\lambda(\mathcal{F}) + \nu(\mathcal{F}) + 3) \cdot \binom{\lambda(\mathcal{F}) + \nu(\mathcal{F})}{\lambda(\mathcal{F})}^2.
\]
\end{lemma}
It follows from~\cite[Lemma 4]{2025Nonpierc} that any non-piercing family \( \mathcal{F} \) of regions satisfies \( \lambda(\mathcal{F}) \le 5 \). Combined with Lemma~\ref{1994}, this gives the bound \( \tau(\mathcal{F}) \le C \nu(\mathcal{F})^{11} \) for some absolute constant \( C > 0 \). Our main result sharpens this by determining the exact maximum value of \( \lambda(\mathcal{F}) \) over all non-piercing families \( \mathcal{F} \).

\begin{figure}[htpb]
    \centering
    \begin{tikzpicture}[scale=1, every node/.style={font=\scriptsize}]
\draw[thick, rotate around={60:(0,2)}] (-1.3,3) ellipse (3 and 0.6);
\node at (-1.5,2) {\(F_1\)};

\draw[thick, rotate around={-60:(0,2)}] (1.3,3) ellipse (3 and 0.6);
\node at (1.5,2) {\(F_2\)};

\draw[thick] (0,-1) ellipse (3 and 0.5);
\node at (0,-1.2) {\(F_3\)};

\draw[thick] (0,0.6) circle (1.5cm);
\node at (0,0.6) {\(F_4\)};
   
\end{tikzpicture}
    \caption{A non-piercing family $\mathcal{F}=\{F_1,F_2,F_3,F_4\}$ of regions with $\lambda(\mathcal{F})=4$.}
    \label{fig:lower bound}
\end{figure}
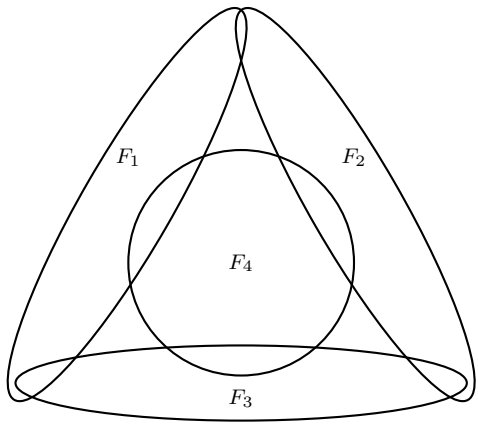

\begin{theorem}\label{thm:lambda=4}
Let \( \mathcal{F} \) be a non-piercing family of regions, then \(\lambda(\mathcal{F}) \le 4.\)
\end{theorem}
\cref{thm:lambda=4}, when combined with~\cref{1994}, directly implies~\cref{thm:main}. Our proof employs a mixture of combinatorial and geometric arguments. Moreover, this bound is tight: there exists a non-piercing family \( \mathcal{F} \) with \( \lambda(\mathcal{F}) = 4 \), see~\cref{fig:lower bound}.

\section{Proof of~\cref{thm:lambda=4}}
\subsection{Topological tools}
First we introduce two classical topological tools that play a central role in our proof: the Jordan Curve Theorem and the strong Hanani–Tutte Theorem. The former states that every simple closed curve in the plane divides the plane into exactly two connected components: one bounded (the interior) and one unbounded (the exterior). A standard consequence is the following:
\begin{lemma}[\cite{Jordan1887}]\label{lem:Jordan}
If two simple closed curves in the plane intersect, then the number of intersection points is even, unless the two curves coincide entirely.
\end{lemma}
The second key tool is a strong variant of the Hanani–Tutte theorem. Informally, it states that every non-planar graph must contain a pair of non-adjacent edges that cross an odd number of times in any drawing. We present the version due to Tutte~\cite{tutte1970towar} as follows:
\begin{theorem}[\cite{tutte1970towar}]\label{tutte}
   In any planar drawing of a subdivision \( G \) of \( K_5 \) or \( K_{3,3} \), there exist two edges, derived from nonadjacent edges of the original \( K_5 \) or \( K_{3,3} \), whose crossing number is odd. 
\end{theorem}

\subsection{High-level overview of the proof}
Although the proof is short, we believe that a brief overview will help highlight its key ideas. As noted earlier, it is straightforward to verify that any non-piercing family \( \mathcal{F} \) satisfies \( \lambda(\mathcal{F}) \le 5 \) based on the proof in~\cite[Lemma 4]{2025Nonpierc}: if every pairwise intersection \( F_i \cap F_j \) contains a point excluded from all other sets, then a simple counting argument reveals a forbidden configuration. However, improving the bound to \( \lambda(\mathcal{F}) \le 4 \) is considerably more subtle. The main challenge is that such configurations are locally indistinguishable from the 5-set case and thus require a more global topological analysis. Assuming \( \lambda(\mathcal{F}) = 5 \), we obtain five sets \( F_1, F_2, F_3, F_4, F_5 \), together with points \(a_{i,j} \in (F_i \cap F_j) \setminus \big(\bigcup_{k \notin \{i,j\}} F_k \big) \quad \text{for all } 1\le i < j\le 5.\)
We take the five points \( a_{i,i+1} \) (indices modulo 5) as the vertices of a hypothetical drawing of the complete graph \( K_5 \), and aim to construct ten simple curves between them to represent the edges of \( K_5 \), each contained in the union of the relevant two sets and avoiding all others.

The core of the argument lies in the careful construction of the connecting curves. For each pair of points \( (a_{i,i+1}, a_{i+1,i+2}) \), we directly select a simple curve in \( F_{i+1} \setminus F_{i+4} \) joining them. In contrast, for each pair \( (a_{i,i+1}, a_{i,i+2}) \), the path must be constructed more delicately: it consists of two concatenated simple curves, lying in \( F_i \setminus F_{i+4} \) and \( F_{i+2} \setminus F_{i+1} \), respectively. The existence of such paths follows from the connectedness of \( F_i \setminus F_j \), a property ensured by the non-piercing condition. To avoid conflicts among the constructed paths, we employ a rotational selection scheme.

This construction yields a drawing of \( K_5 \) in which many pairs of independent edges are guaranteed not to cross. However, by the strong Hanani--Tutte theorem, such a drawing must contain a pair of independent edges that cross an odd number of times. By identifying such a pair and showing that they can be extended to form two Jordan curves intersecting an odd number of times, we arrive at a contradiction with the Jordan curve theorem.

The novelty of the argument lies in the precise geometric construction of the curves, which leverages the non-piercing property to carefully manage their placement and intersection behavior. This approach bridges local structural constraints with global topological tools and is essential to rule out the existence of such five sets.
\subsection{Formal Proof}
We proceed by contradiction. Throughout the proof, all indices are considered modulo $5$. Suppose, for contradiction, that there exist five sets \( F_1, F_2, \ldots, F_5 \in \mathcal{F} \) such that for every pair \( \{i,j\} \in \binom{[5]}{2} \), there is a distinguished point \( a_{i,j} \in F_i \cap F_j \) with \( a_{i,j} \notin F_k \) for any \( k \in [5] \setminus \{i,j\} \). We then select some of these points to define the five vertices of our \( K_5 \) construction. For each \( i \in [5] \), let \( v_i = a_{i,i+1} \). By definition, \( v_i \in F_i \cap F_{i+1} \) and \( v_i \notin F_{i+2} \cup F_{i+3} \cup F_{i+4} \).

Next, for each pair \( \{i,j\} \in \binom{[5]}{2} \), we define a simple curve \( p_{i,j} \) connecting \( v_i \) and \( v_j \) as follows. Since indices are taken modulo 5, each such pair satisfies either \( j \equiv i+1 \pmod{5} \) or \( j \equiv i+2 \pmod{5} \) (up to swapping $i$ and $j$). Thus, it suffices to define the curves in these two cases:

\begin{itemize}
    \item For each \( i \in [5] \), define \( p_{i,i+1} \) as a simple curve connecting \( v_i \) to \( v_{i+1} \) entirely contained in \( F_{i+1} \setminus F_{i+4} \). Such a curve exists since both endpoints lie in the region $F_{i+1} \setminus F_{i+4}$, and $F_{i+1}\setminus F_{i+4}$ is connected by non-piercing property.
    
    \item For each \( i \in [5] \), define \( p_{i,i+2} \) as a concatenation of two simple curves passing through the point \( a_{i,i+2} \): one from \( v_i \) to \( a_{i,i+2} \) inside \( F_i \setminus F_{i+4} \), and another from \( a_{i,i+2} \) to \( v_{i+2} \) inside \( F_{i+2} \setminus F_{i+1} \). The existence of these curves follows from $v_i, a_{i,i+2} \in F_{i} \setminus F_{i+4}$, $a_{i,i+2}, v_{i+2} \in F_{i+2} \setminus F_{i+1}$, and $F_{i}\setminus F_{i+4}$, $F_{i+2}\setminus F_{i+1}$ are both connected by non-piercing property.
\end{itemize}

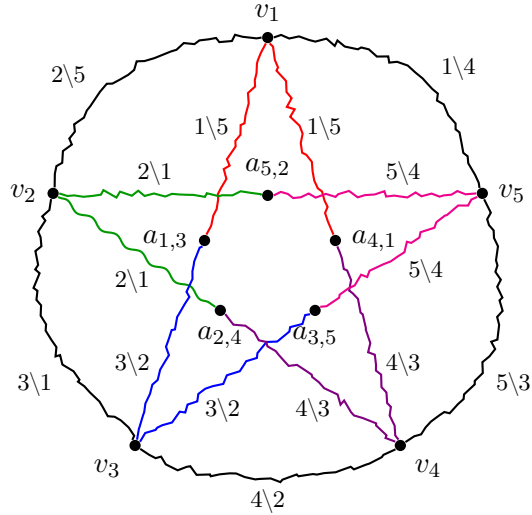
\begin{figure}[htpb]
    \centering
   \begin{tikzpicture}[scale=3, every node/.style={font=\normalsize}]

\node[fill=black, circle, inner sep=1.5pt, label=above:$v_1$] (a12) at (90:1) {};
\node[fill=black, circle, inner sep=1.5pt, label=left:$v_2$] (a23) at (162:1) {};
\node[fill=black, circle, inner sep=1.5pt, label=below left:$v_3$] (a34) at (234:1) {};
\node[fill=black, circle, inner sep=1.5pt, label=below right:$v_4$] (a45) at (306:1) {};
\node[fill=black, circle, inner sep=1.5pt, label=right:$v_5$] (a51) at (18:1) {};

\node[fill=black, circle, inner sep=1.5pt, label=above:$a_{5,2}$] (a13) at (0,0.3) {};
\node[fill=black, circle, inner sep=1.5pt, label=left:$a_{1,3}$] (a24) at (-0.28,0.1) {};
\node[fill=black, circle, inner sep=1.5pt, label=below:$a_{2,4}$] (a35) at (-0.21,-0.21) {};
\node[fill=black, circle, inner sep=1.5pt, label=below:$a_{3,5}$] (a41) at (0.21,-0.21) {};
\node[fill=black, circle, inner sep=1.5pt, label=right:$a_{4,1}$] (a52) at (0.3,0.1) {};

 
\draw[thick, decorate, decoration={random steps, segment length=3pt, amplitude=1.5pt}]
  (a12) .. controls (-0.8,1) and (-0.9,0.5).. (a23)
  node[midway, above left, font=\footnotesize] {$2\setminus 5$};

\draw[thick, decorate, decoration={random steps, segment length=3pt, amplitude=1.5pt}] 
(a23) .. controls (-1.2,-0.4) and (-0.7,-0.6).. (a34) 
node[midway, below left, font=\footnotesize] {$3\setminus 1$};

\draw[thick, decorate, decoration={random steps, segment length=3pt, amplitude=1.5pt}]
(a34) .. controls (-0.3,-1) and (0.3,-1).. (a45) 
node[midway, below, font=\footnotesize] {$4\setminus 2$};

\draw[thick, decorate, decoration={random steps, segment length=3pt, amplitude=1.5pt}]
(a45) .. controls (0.9,-0.6) and (1.15,-0.4).. (a51) 
node[midway, below right, font=\footnotesize] {$5\setminus 3$};

\draw[thick, decorate, decoration={random steps, segment length=3pt, amplitude=1.5pt}]
(a51) .. controls (0.8,0.6) and (0.8,1).. (a12) node[midway, above right, font=\footnotesize] {$1\setminus 4$};


\draw[thick, decorate, decoration={random steps, segment length=3pt, amplitude=1.5pt}, red]
 (a12) -- (a24);
\node at (-0.25,0.6) {\footnotesize $1\setminus 5$};

\draw[thick, decorate, decoration={random steps, segment length=3pt, amplitude=1.5pt}, red] 
(a12) -- (a52);
\node at (0.25,0.6) {\footnotesize $1\setminus 5$};

\draw[thick, decorate, decoration={snake, amplitude=0.35mm, segment length=4mm}, green!60!black] (a23) -- (a35);
\node at (-0.6,-0.08) {\footnotesize $2\setminus 1$};

\draw[thick, decorate, decoration={random steps, segment length=3pt, amplitude=1.5pt}, green!60!black] (a23) -- (a13);
\node at (-0.5,0.4) {\footnotesize $2\setminus 1$};

\draw[thick, decorate, decoration={random steps, segment length=3pt, amplitude=1.5pt}, blue] (a34) -- (a24);
\node at (-0.6,-0.45) {\footnotesize $3\setminus 2$};

\draw[thick, decorate, decoration={random steps, segment length=3pt, amplitude=1.5pt}, blue] (a34) -- (a41);
\node at (-0.2,-0.65) {\footnotesize $3\setminus 2$};

\draw[thick, decorate, decoration={random steps, segment length=3pt, amplitude=1.5pt}, violet] (a45) -- (a52);
\node at (0.6,-0.45) {\footnotesize $4\setminus 3$};

\draw[thick, decorate, decoration={random steps, segment length=3pt, amplitude=1.5pt},violet] (a35) -- (a45);
\node at (0.2,-0.65) {\footnotesize $4\setminus 3$};

\draw[thick, decorate, decoration={random steps, segment length=3pt, amplitude=1.5pt}, magenta] (a51) -- (a41);
\node at (0.7,-0.035) {\footnotesize $5\setminus 4$};

\draw[thick, decorate, decoration={random steps, segment length=3pt, amplitude=1.5pt}, magenta] (a51) -- (a13);
\node at (0.6,0.4) {\footnotesize $5\setminus 4$};

\end{tikzpicture}
    \caption{A planar embedding of the complete graph $K_5$ with vertices $\{v_1,v_2,v_3,v_4,v_5\}$.}
    \label{fig:upper bound}
\end{figure}

Note that the pair \( \{i,j\} \) is unordered, so \( a_{i,j} = a_{j,i} \) and \( p_{i,j} = p_{j,i} \). We illustrate our embedding of \( K_5 \) in~\cref{fig:upper bound}, where the notation \( a \setminus b \) indicates a simple curve contained in \( F_a \setminus F_b \) for \( a,b \in [5] \).

\begin{claim}\label{claim:Verification}
    For any four distinct indices $\{i, j, k, \ell\} \subseteq [5]$, the crossing number of $p_{i,j}$ and $p_{k,\ell}$ is even.
\end{claim}
\begin{poc}
By symmetry, it suffices to consider $(i,j,k,\ell) = (1,2,3,4)$, $(1,3,2,4)$, and $(1,4,2,3)$.
\begin{itemize}
    \item 
For $(i,j,k,\ell) = (1,2,3,4)$: We have $p_{1,2} \subseteq F_2 \setminus F_5$ and $p_{3,4} \subseteq F_4 \setminus F_2$, thus they are disjoint and their crossing number is $0$.

    \item 
For $(i,j,k,\ell) = (1,4,2,3)$: The curve $p_{2,3}$ lies entirely in $F_3 \setminus F_1$. The curve $p_{1,4}$ consists of two simple curves: one in $F_4 \setminus F_3$ and another in $F_1 \setminus F_5$. Since both of these sets are disjoint from $F_3 \setminus F_1$, neither simple curve in $p_{1,4}$ intersects $p_{2,3}$. Therefore, their crossing number is $0$.

    \item 
For $(i,j,k,\ell) = (1,3,2,4)$: The curve $p_{1,3}$ consists of two simple curves, one in $F_1 \setminus F_5$ and another in $F_3 \setminus F_2$. The curve $p_{2,4}$ consists of two simple curves, one in $F_2 \setminus F_1$ and another in $F_4 \setminus F_3$. Since $F_2 \setminus F_1$ is disjoint from both $F_1 \setminus F_5$ and $F_3 \setminus F_2$, and $F_3 \setminus F_2$ is disjoint from $F_4 \setminus F_3$, the only possible intersection occurs between the simple curve in $F_1 \setminus F_5$ and the one in $F_4 \setminus F_3$.

Let $\alpha$ be the simple curve in $F_4 \setminus F_3$ from $v_4$ to $a_{2,4}$, and let $\beta$ be the simple curve in $F_1 \setminus F_5$ from $v_1$ to $a_{1,3}$. 
Let $\alpha'$ be a simple curve in $F_4 \setminus F_1$ from $v_4$ to $a_{2,4}$, and $\beta'$ a simple curve in $F_1 \setminus F_4$ from $v_1$ to $a_{1,3}$, their existence being guaranteed because of the non-piercing property since $v_4,a_{2,4} \in F_4 \setminus F_1$ and $v_1,a_{1,3} \in F_1 \setminus F_4$.
Also notice that $\alpha \cup \alpha'$ and $\beta \cup \beta'$ form simple closed curves. Since $\alpha,\alpha' \subseteq F_4$, $\beta,\beta' \subseteq F_1$, and $\alpha' \cap F_1 = \emptyset = \beta' \cap F_4$, it follows that the intersection number of $\alpha$ and $\beta$  equals that of $\alpha \cup \alpha'$ and $\beta \cup \beta'$, which is even by~\cref{lem:Jordan}.
\end{itemize}\end{poc}
However, by~\cref{tutte},~\cref{claim:Verification} cannot occur. This finishes the proof.

\section{Concluding remarks}
The result of~\cref{1994} has found numerous applications in the literature, including the influential work of Fox, Pach, and Suk on the sunflower conjecture for set systems with bounded VC-dimension~\cite{2023BVCSunflower}, as well as the resolution of Scott’s induced subdivision conjecture for maximal triangle-free graphs~\cite{2012Nicolas}. However, to the best of our knowledge, exact values for the largest dyadic dual VC-dimension of specific geometric set systems, such as the one established in this paper, are rare. This suggests that determining the dyadic dual VC-dimension of a given set system may be an interesting and worthwhile goal in its own right, especially for set systems arsing from geometry.

In particular, the dyadic dual VC-dimension $\lambda(\mathcal{F})$ is often difficult to control. For instance, if a set system has VC-dimension at most 1, then as observed in~\cite{2023BVCSunflower}, it is known that $\lambda(\mathcal{F}) \leq 3$. Yet, once the VC-dimension exceeds 1, no general bound on $\lambda(\mathcal{F})$ is known, even when the VC-dimension is as low as 2. This reflects a striking gap between the classical VC-dimension and its dyadic dual counterpart, and underscores the subtlety of this parameter.

Another direction that we find particularly intriguing is whether the quantitative bound in~\cref{1994} can be further improved. Even modest improvements could have meaningful consequences for transversal-type problems in geometric and combinatorial settings.

\section*{Acknowledgement}
Xinqi Huang was supported by the National Key Research and Development Programs of China 2023YFA1010200 and 2020YFA0713100, the NSFC under Grants No. 12171452 and No. 12231014, the Innovation Program for Quantum Science and Technology 2021ZD0302902, and the Institute for Basic Science (IBS-R029-C4). 
Yuzhen Qi was supported by China Scholarship Council and IBS-R029-C4. 
Mingyuan Rong was supported by National Key Research and Development Program of China 2023YFA1010201, the NSFC under Grants No. 12125106 and the Excellent PhD Students Overseas Study Program of the University of Science and Technology of China.
Zixiang Xu was supported by the Institute for Basic Science (IBS-R029-C4) and would like to thank Minho Cho for helpful discussion.
\bibliographystyle{abbrv}
\bibliography{Nonpiercing}
\end{document}